\documentclass[12pt,a4paper]{article}

\newcommand{\bQ}{\mathbb{Q}}
\newcommand{\bR}{\mathbb{R}}
\newcommand{\bC}{\mathbb{C}}
\newcommand{\bN}{\mathbb{N}}
\newcommand{\bD}{\mathbb{D}} 
\newcommand{\hbC}{\hat{\mathbb{C}}}
\newcommand{\supp}{\operatorname{supp}}

\usepackage{amsmath,amsthm,amssymb}

\theoremstyle{plain}
\newtheorem{problem}{Problem}

\newtheorem{theorem}{Theorem}[section]

\newtheorem{lemma}{Lemma}[section]
\newtheorem{proposition}{Proposition}[section]
\newtheorem{mainth}{Main Theorem}
\theoremstyle{remark}
\newtheorem{remark}{Remark}

\newtheorem{pfcritical}{Proof of Main Theorem $\ref{mainth:critical}$}
 
\newtheorem{pfsingularity}{Proof of Main Theorem $\ref{mainth:singularity}$}

\newtheorem{pfpreserve}{Proof of Proposition $\ref{prop:preserve}$}

\theoremstyle{definition}
\newtheorem{definition}{Definition}[section]
\newtheorem{C}{Convention}

\newtheorem{ac}{ACKNOWLEDGMENT}

\begin{document} 

\title{Linearization problem on structurally finite entire functions} 

\author{Y\^usuke Okuyama\footnote{Partially supported by the Sumitomo
Foundation, and the Ministry of Education, Science, Sports and Culture,
Grant-in-Aid for Young Scientists (B), 15740085, 2004.}\\
Department of Mathematics, Faculty of Science, \\
Kanazawa University, Kanazawa 920-1192 Japan\\
email; okuyama@kenroku.kanazawa-u.ac.jp\\
\ \\
\it Dedicated to the memory of Professor Nobuyuki Suita}



\date{}

\maketitle

\begin{abstract}
 We show that if a 1-hyperbolic structurally finite entire function
 of type $(p,q)$, $p\ge 1$, is linearizable at an irrationally
 indifferent fixed point, then its multiplier satisfies the Brjuno
 condition. We also prove the generalized Ma\~n\'e theorem;
 if an entire function has only finitely many critical points and
 asymptotic values, then for every such a non-expanding forward invariant set
 that is either a Cremer cycle or the boundary of a cycle of Siegel disks,
 there exists an asymptotic value or a recurrent critical point such that
 the derived set of its forward orbit contains this invariant set.
 From it, the concept of $n$-subhyperbolicity naturally arises.
\end{abstract}

\section{Introduction}\label{sec:main}
A structurally finite entire function is constructed from 
finitely many quadratic blocks and exponential blocks by 
{\it Maskit surgeries} which connect two functions.

\begin{definition}[structural finiteness]
 A {\it structurally finite entire function of type} $(p,q)$ is an entire
 function constructed from $p$ {\itshape quadratic polynomial blocks}
 and $q$ {\itshape exponential function blocks}. $\mathcal{SF}_{p,q}$
 denotes the set of all structurally finite entire functions of type $(p,q)$.
\end{definition}

For the precise definition and details,
see \cite{tanighongkong} or \cite{tanig01}.
By this definition, we have that
a structurally finite entire function has in fact
both the {\itshape topological
characterization} and the {\itshape explicit representation}, which 
are used in Section \ref{sec:perturb} and \ref{sec:proof}
respectively:
\begin{theorem}[topological characterization \cite{tanighongkong}]\label{th:topological}
 Every element of $\mathcal{SF}_{p,q}$
 has exactly $p$ critical points and $q$ transcendental singularities of
 its inverse.
 Conversely, every entire function with exactly $p$ critical points and
 $q$ transcendental singularities of its inverse
 belongs to $\mathcal{SF}_{p,q}$.
\end{theorem}
In the above and the whole paper,
we always count critical points {\itshape with} multiplicities.
Hence structurally finite entire functions are of the Speiser class,
that is, have finitely many singular values.
The classification theorem of Fatou components of
this class is known (\cite{GK86} and \cite{EL92}).
In particular, there are neither wandering nor Baker domains.

\begin{theorem}[explicit representation \cite{tanig01}]\label{th:explicit}
For $(p,q)\neq(0,0)$, $\mathcal{SF}_{p,q}$
agrees with $SF_{p,q}$, where
\[
 SF_{p,q}:=\left\{\int_0^z(c_pt^p+\cdots+c_0)e^{a_qt^q+\cdots+a_1t}dt+b;
c_pa_q\neq 0\right\}\quad (q\neq 0),
\]
and $SF_{p,0}:={\rm Poly}_{p+1}=\{\text{polynomials of degree }p+1\}$.
\end{theorem}

From now on, we assume $\lambda=e^{2\pi i\alpha}$ $(\alpha\in\bR-\bQ)$.

Let us consider an {\itshape irrationally indifferent} cycle
of an entire function $f$ of period $n$ with multiplier $\lambda$.
It is called a {\it Siegel} cycle
if every point of this cycle has a neighborhood where
the {\itshape first return map} $f^n$ is conformally conjugate to
$R_{\lambda}(z)=\lambda z$ on the unit disk.
Otherwise it is called a {\it Cremer} cycle.

The Brjuno condition for $\alpha$ means that
\[
 \sum_{n=0}^{\infty}\frac{\log{q_{n+1}}}{q_n}<\infty,
\]
where $\{p_n/q_n\}$ is the sequence of
rational numbers approximating $\alpha$ defined
by its continued fraction expansion.

The following shows that an irrationally indifferent cycle
is Siegel when $\alpha$ satisfies this Brjuno condition:
\begin{theorem}[Brjuno \cite{Brjuno}]
 Let $f(z)=\lambda z+\cdots$ be an analytic germ at the origin.
 If $\alpha$ satisfies the Brjuno condition,
 then $f$ is $($analytically$)$ linearizable, that is,
 on a neighborhood of the origin, $f$ is conformally conjugate to
 $R_{\lambda}(z)=\lambda z$ on the unit disk.
\end{theorem}

In \cite{Yoccoz96},
Yoccoz gave a beautiful alternative proof of this Brjuno theorem
and also showed the following theorem in the case of period one.
Later we generalized it in the case of arbitrary period:
\begin{theorem}[Yoccoz \cite{Yoccoz96}, Okuyama \cite{1-subhyp}]\label{th:quad}
 If an irrationally indifferent cycle of a quadratic polynomial 
 with multiplier $\lambda$ is a Siegel cycle,
 then $\alpha$ satisfies the Brjuno condition.
\end{theorem}

Even in the cubic polynomial case,
it is not known whether Theorem \ref{th:quad} can be generalized.
In the transcendental entire function case, Lukas Geyer showed
the following:

\begin{theorem}[Geyer \cite{Geyer01}]
 If the origin is a Siegel fixed point of $\lambda\int^z_0(1+t)e^tdt$,
 then $\alpha$ satisfies the Brjuno condition.
\end{theorem}

Clearly Geyer's example belongs to $\mathcal{SF}_{1,1}$,
quadratic polynomials to $\mathcal{SF}_{1,0}$, and both of them
naturally satisfy the $1$-{\itshape hyperbolicity} defined in
Section \ref{sec:subhyp} (see also \cite{1-subhyp}).
In this paper, we shall extend Geyer's result
to $1$-{\itshape hyperbolic} structurally finite entire functions 
by more general and synthetic method:

\begin{mainth}\label{mainth:critical}
 If a $1$-hyperbolic structurally finite entire function of type
 $(p,q)$, $p\ge 1$, has a Siegel fixed point with multiplier
 $\lambda=e^{2\pi i\alpha}$, then $\alpha$ satisfies the Brjuno condition.   
\end{mainth}

In the case $p=0$, we just have:
\begin{mainth}\label{mainth:singularity}
 If a $1$-hyperbolic structurally finite entire function of type
 $(0,q)$, $q\ge 1$, has a Siegel fixed point with multiplier
 $\lambda=e^{2\pi i\alpha}$, then $E(z)=\lambda\int_0^ze^tdt$ 
 is linearizable at the origin.
\end{mainth}

We note that $E(z)=\lambda\int_0^ze^tdt\in\mathcal{SF}_{0,1}$.
Hence the most fundamental case remains open:
\begin{problem}
 If the origin is a Siegel fixed point of $E(z)=\lambda\int_0^ze^tdt$,
 then $\alpha$ satisfies the Brjuno condition?
\end{problem}

In Section \ref{sec:subhyp}, we define the {\it $n$-subhyperbolicity}
in the way similar to that in \cite{1-subhyp}.
The {\itshape generalized} Ma\~n\'e theorem is crucial.
In Section \ref{sec:perturb}, we shall explain the {\it 
linearizability-preserving perturbation} of the $n$-hyperbolic
entire function, which increases the number of the foliated
equivalence classes of acyclic singular values in the Fatou set.
The topological characterization of $\mathcal{SF}_{p,q}$
shows that this perturbation is {\itshape closed} in
$\mathcal{SF}_{p,q}$.
In Section \ref{sec:proof}, we shall prove Main Theorems.
Using the explicit representation of $\mathcal{SF}_{p,q}$,
we can apply the {\itshape algebraic quadratic perturbation},
which was first applied to polynomials by P\'erez-Marco \cite{PMar},
to structurally finite entire functions.
In Section \ref{sec:mane}, we shall give a proof of the generalized 
Ma\~n\'e theorem, from which we naturally derive
the concept of $n$-subhyperbolicity.

\begin{ac}
 The author would like to express his gratitude to
 Prof.\ Masahiko Taniguchi and Prof.\ Toshiyuki Sugawa
 for valuable discussions and advices, and
 to Prof.\ Mitsuhiro Shishikura for commenting
 on the proof of the corollary of the Ma\~n\'e theorem.
 
 The author must particularly thanks the referee. The referee read the paper
 with great care, pointed out many minor mistakes, the important one of which 
 was involved with the generalized Ma\~n\'e theorem; 
 the author could correct the statement of Theorem \ref{th:point}
 and refine the proof of Theorem \ref{th:correspond}. 
 The referee's careful comments also helped
 the author refine the presentation of this article.
\end{ac}

\section{$n$-subhyperbolicity}\label{sec:subhyp}
We assume that entire functions are neither constant nor linear.
Let $f$ be an entire function and $F(f)$ and $J(f)$ 
the Fatou and Julia sets of $f$ respectively. 

\begin{definition}[derived set, $\omega$-limit set, and recurrence]
The {\it derived set} $d(c)$ of $c\in\bC$ is defined by the set of
all derived (or accumulation) points $z\in\bC$ of $\{f^n(c)\}_{n\in\bN}$,
i.e., $z\in\bC$ such that for every neighborhood $U$ of $z$,
$(U-\{z\})\cap\{f^n(c)\}_{n\in\bN}\neq\emptyset$.
The {\it $\omega$-limit set} $\omega(c)$ of $c\in\bC$ is defined 
by the set of all $z\in\bC$ such that $\lim_{i\to\infty}f^{n_i}(c)=z$
for some increasing $\{n_i\}\subset\bN$.
A point $c$ is {\it recurrent} if $\omega(c)\ni c$. 
\end{definition}
\begin{remark}
 If $c$ is either periodic or preperiodic,
 then $d(c)=\emptyset$ and $\omega(c)$ equals
 the cycle where $c$ is eventually mapped. 
 Otherwise $d(c)=\omega(c)$.
 In particular, if $c$ is a recurrent critical
 point of $f$ in $J(f)$, $d(c)=\omega(c)$.
\end{remark}

\begin{definition}[correspondence]
An asymptotic value or a recurrent critical point $s$
{\it corresponds to} an irrationally indifferent cycle $C$
if $d(s)\supset\Gamma$. Here 
let $\Gamma=\Gamma_C$ be,  if $C$ is Siegel, the boundary of the cycle of
the Siegel disks associated with $C$, otherwise the cycle $C$ itself. 
\end{definition}

\begin{remark}
 If $s$ corresponds to an irrationally indifferent cycle, then
 $d(s)\neq\emptyset$, so $s$ is neither periodic nor preperiodic.
 It also holds that $s\in J(f)$.
\end{remark}

We shall show the following in Section \ref{sec:mane}. 
\begin{theorem}[the generalized Ma\~n\'e theorem]\label{th:correspond}
 Suppose that $f$ has only finitely many critical points and asymptotic
 values. Then for every irrationally indifferent cycle $C$,
 there exists an asymptotic value or a recurrent
 critical point corresponding to $C$.
\end{theorem}

We fix the definition of the {\it transcendental singularities}
of the inverse of an entire function $f$.
For $a\in\bC$, let $A:=\{A(r)\}_{r>0}$ be a family of domains in $\bC$
such that for $r>0$, $A(r)$ is a component of $f^{-1}(\bD_r(a))$
and if $0<r_1<r_2$, then $A(r_1)\subset A(r_2)$.
Then the intersection of all the closures of $A(r)$ in $\hat{\bC}$
consists of only one point. If this point is the infinity, $A$ is
called a {\it transcendental singularity of} $f^{-1}$ over $a$.
We note that then the $a$ is an asymptotic value of $f$, and that
the number of transcendental singularities of $f^{-1}$
is not less than that of asymptotic values.

\begin{definition}[correspondence]
The transcendental singularity $A$ of $f^{-1}$ over the asymptotic value
$a$ {\it corresponds to} $\Gamma$ if so does $a$. 
\end{definition}

\begin{C}
 For a transcendental singularity $A$ of $f^{-1}$ over $a$, 
 we say that the image of $A$ by $f$ is $a$, and write $f(A)=a$.
 Moreover, if $a\in J(f)$, we say that $A\in J(f)$.
\end{C}

Now we define the {\it $n$-subhyperbolicity}.

\begin{definition}[$n$-subhyperbolicity]
Let $f$ be a structurally finite entire function.
For a non-negative integer $n$, $f$ is {\it $n$-subhyperbolic} if 
\begin{enumerate}
\def\labelenumi{(\roman{enumi})}
\item there exist exactly $n$ recurrent critical points of $f$ or
      transcendental singularities of $f^{-1}$ each of which corresponds to
      some irrationally indifferent cycle of $f$,
\item every critical point of $f$ and transcendental singularity of $f^{-1}$
      in $J(f)$ other than such ones as in (i) is preperiodic, and 
\item no orbits of singular values in $F(f)$ accumulate to $J(f)$.
\end{enumerate}
 An $n$-subhyperbolic $f$ is {\it $n$-hyperbolic} if it has no such
 (ideal) points as in (ii). 
\end{definition}

In the case $f$ is a polynomial, this definition agrees with
what we have defined in \cite{1-subhyp}.
For several examples of $n$-subhyperbolic
polynomials, see \cite{1-subhyp}.

\section{Linearizability preserving perturbation}\label{sec:perturb}
From now on, we fix $\lambda=e^{2\pi i\alpha}$, where $\alpha\in\bR-\bQ$.

For an entire function $f$, 
a point is said to be {\it acyclic} if it is neither periodic nor
preperiodic point of $f$. The {\it grand orbit} of $x\in\bC$ is
the set
\[
 \{y\in \bC;f^i(x)=f^j(y)\text{ for some $i,j\ge 0$}\}.
\]
$x,y\in\bC$ are in the {\it foliated equivalence class} of
$f$ if the closure of their grand orbits agree with each other.
Let $N_{AC}(f)$ be the number of the foliated equivalence classes of 
acyclic singular values of $f$ in $F(f)$. 

\begin{proposition}[linearizability preserving perturbation]\label{prop:preserve}
 Let $f\in\mathcal{SF}_{p,q}$ be $n$-hyperbolic and have 
 a Siegel fixed point with multiplier $\lambda$.
 Then there exists an $n$-hyperbolic $g\in\mathcal{SF}_{p,q}$ such that
\begin{enumerate}
\def\labelenumi{$($\roman{enumi}$)$}
 \item $g$ also has a Siegel fixed point with multiplier $\lambda$, and
 \item $N_{AC}(g)=p+q-n$.
\end{enumerate}
\end{proposition}

In the rest of this section, we prove Proposition \ref{prop:preserve}.
Let $f\in\mathcal{SF}_{p,q}$.
We use the following lemmas,
the first three of which are about perturbations of critical points,
and essentially proved in Section 2 in \cite{1-subhyp}.
See also Section 5 in \cite{tanighongkong}, ``relaxing the relations
between the singularity data''. 

\begin{lemma}\label{lem:decompose}
 Let $c$ be a non-periodic critical point in $F(f)$ with multiplicity
 $m\ge 2$. 
 There exists a Jordan neighborhood $U$ of $c$ in $F(f)$
 such that $\overline{U}-\{c\}$ contains no critical point,
 $f$ maps $U$ onto some Jordan domain properly,
 and $U\cap\bigcup_{n\ge 1}f^{n}(U)=\emptyset$.
 And there exist a quasiconformal automorphism $\Phi$ of
 $\bC$ and $g\in\mathcal{SF}_{p,q}$ such that 
 $g$ has exactly $m$ distinct critical points in
 $\Phi(U)$, which are simple, and 
 $g=\Phi\circ f\circ\Phi^{-1}$ on $\bC- \Phi(U)$.
\end{lemma}

We note that the following 
Lemma \ref{lem:decomposeperiodic} is the inverting Carleson
and Gamelin operation in \cite{CG}.
\begin{lemma}\label{lem:decomposeperiodic}
 Let $c$ be a periodic critical point in $F(f)$ with multiplicity $m\ge 1$ and  of period $p$.
 There exists a Jordan neighborhood $U$ of $c$ in $F(f)$
 such that $\overline{U}$ contains no critical point of $f$ other than $c$,
 $f$ maps $U$ onto some Jordan domain properly,
 and $f^{p}(U)\Subset U$.
 And there exist a quasiconformal automorphism $\Phi$ of
 $\mathbb{C}$ and $g\in\mathcal{SF}_{p,q}$ such that
 $g$ has exactly $m$ distinct critical points
 in $\Phi(U)$, which are simple, $\Phi(c)$ is not a critical point of
 $g$, and $g=\Phi\circ f\circ\Phi^{-1}$ on $\bC-\Phi(U-\{c\})$.
\end{lemma}

\begin{lemma}\label{lem:move}
 Let $c$ be a non-periodic and simple critical point in $F(f)$.
 There exists a Jordan neighborhood $U$ of $c$ in $F(f)$ such that
 $\overline{U}-\{c\}$ contains no critical point,
 $f$ maps $U$ onto some Jordan domain properly,
 and $U\cap\bigcup_{n\ge 1}f^{n}(U)=\emptyset$.
 For every $y\in f(U)$,
 there exist a quasiconformal automorphism $\Phi$ of $\bC$
 and $g\in\mathcal{SF}_{p,q}$ such that
 $g$ has only one critical point $\Phi(c)$ in $\Phi(U)$,
 which is simple, $g(\Phi(c))=\Phi(y)$, and
 $g=\Phi\circ f\circ\Phi^{-1}$ on $\bC-\Phi(U)$.
\end{lemma}

We need to show the following here. 
In the proof, the classification theorem of Fatou components
already stated in Section \ref{sec:main} is implicitly used.
\begin{lemma}\label{lem:transmove}
 Let $A=\{A(r)\}_{r>0}$ be a transcendental singularity of $f^{-1}$
 over $a\in F(f)$. There exists an $r_1>0$ such that
 $U:=A(r_1)$ is in $F(f)$, disjoint from
 $\bigcup_{n\ge 0}f^n(\bD_{r_1}(a))$, and $\overline{U}$
 contains no critical points of $f$.
 For every $\epsilon\in\bD_{r_1/2}$, 
 there exist quasiconformal automorphisms $\Psi$ and $\Phi$ of $\bC$
 and $g\in\mathcal{SF}_{p,q}$ such that
 $\Psi$ is the identity outside $\bD_{r_1}(a)$, $\Psi(a)=a+\epsilon$,
 and $g=\Phi\circ\Psi\circ f\circ\Phi^{-1}$ on $\bC$. In particular,
 $g=\Phi\circ f\circ\Phi^{-1}$ on $\bC-\Phi(U)$ and 
 $\Phi(a+\epsilon)$ is an asymptotic value of $g$.
\end{lemma}

\begin{remark}
 In the above Lemmas, 
 we can assume that the diameter of $U$ is arbitrarily small.
\end{remark}

\begin{proof}
 Since $\bigcup_{n\ge 0}f^n(\bD_r(a))$ is uniformly bounded for 
 sufficiently small $r>0$, 
 there exists so small $r_1>0$ that $A(r_1)\cap\bigcup_{n\ge 0}f^n(\bD_{r_1}(a))=\emptyset$ and $\overline{A(r_1)}$ contains no
 critical points. Let $\eta:[0,+\infty)\to[0,1]$ be such a smooth
 function that identically equals one and zero on $[0,r_1/2]$
 and on $[r_1,+\infty)$ respectively.
 For $\epsilon\in\bD_{r_1/2}$, we put $\Psi(w):=w+\epsilon\eta(|w-a|)$.
 It is easy to see that $\tilde{f}:=\Psi\circ f$
 is a quasiregular map on $\bC$,
 equals $f+\epsilon$ and $f$ on $A(r_1/2)$ and on $\bC-A(r_1)$ respectively,
 and is unbranched on $\overline{A(r_1)}-A(r_1/2)$. 
 Let $\mu$ be the Beltrami coefficient on $A(r_1)$ of $\tilde{f}$ and
 we define $\tilde{\mu}$ by the pullback $(f^*)^n\mu$ on
 $f^{-n}(A(r_1))$ ($n\in\bN\cup\{0\}$) and 0 on $\bigcap_{n\ge 0}(\bC-
 f^{-n}(A(r_1)))$, which is an $\tilde{f}$-invariant
 Beltrami coefficient on $\bC$.
 Let $\Phi$ be a quasiconformal automorphism of $\bC$ whose Beltrami
 differential equals $\tilde{\mu}$. Then $g=\Phi\circ\tilde{f}\circ\Phi^{-1}$
 is an entire function. Finally, 
 since $g$ has the same number of critical points and transcendental
 singularities as $f$, it follows that $g$ is in $\mathcal{SF}_{p,q}$
 from Theorem \ref{th:topological}, the topological characterization 
 of structurally finite entire functions.
\end{proof}

\begin{pfpreserve}
 Suppose that $f$ is $n$-hyperbolic and have a Siegel fixed point $z_0$
 with multiplier $\lambda$. 
 By applying the above perturbations to $f$ inductively and 
 in finitely many times,
 we obtain $g\in\mathcal{SF}_{p,q}$, which
 satisfies $N_{AC}(g)=p+q-n$, is $n$-hyperbolic, and 
 is quasiconformally conjugate to $f$ around $z_0$ on a
 neighborhood of $\Phi(z_0)$. 
 Therefore $g$ has a Siegel fixed point $\Phi(z_0)$ with multiplier
 $\lambda$ (see p.~61--p.~62 in \cite{Yoccoz96}). \qed
\end{pfpreserve}

\section{Proof of Main Theorems}\label{sec:proof}
Let $f\in \mathcal{SF}_{p,q}$ have an irrationally indifferent fixed
point $z_0$ with multiplier $\lambda=e^{2\pi i\alpha}$, where
$\alpha\in\bR-\bQ$. In the case $q=0$, Main theorem
\ref{mainth:critical} is proved in \cite{1-subhyp}.
Therefore we assume $q\ge 1$.

By Theorem \ref{th:explicit}, that is, the explicit representation,
and by an affine conjugation which maps $z_0$ to the origin,
we assume that 
\begin{equation}
 f(z)=\lambda\int_0^zP(t)e^{Q(t)}dt, 
\end{equation}
where $P$ is a polynomial of degree $p$ with $P(0)=1$
and $Q$ is that of degree $q$ with $Q(0)=0$. 
Let $SF_{p,q}(\lambda)$ be the set of all such functions,
which is a $(p+q)$-dimensional complex manifold with respect to
coefficients of $P$ and $Q$. Furthermore,
we say $f_1\sim f_2$ ($f_1,f_2\in SF_{p,q}(\lambda)$)
if $f_1(cz)/c=f_2$ for some $c\in\bC^*$.

Suppose that $f$ is 1-hyperbolic and the origin is a Siegel fixed point.
By Proposition \ref{prop:preserve}, we can also assume that
$N_{AC}=p+q-1$, which equals the complex 
dimension of $SF_{p,q}(\lambda)/\sim$.
Since $f$ is 1-hyperbolic, it has no parabolic cycle.
Therefore, by the same argument as that in Lemma 4.1 in 
\cite{1-subhyp} (see also \cite{MS94}), 
the image of the uniformization map (holomorphic injection)
from the Teichm\"uller space of $f$ into $SF_{p,q}(\lambda)/\sim$
becomes a {\itshape domain} in $SF_{p,q}(\lambda)/\sim$.
Hence $f$ is {\it quasiconformally stable}
in $SF_{p,q}(\lambda)$, that is, there exists an open neighborhood
of $f$ every element of which is quasiconformally conjugate to $f$.

\begin{pfcritical}
Suppose that $p\ge 1$. 
Since $f$ is quasiconformally stable in $SF_{p,q}(\lambda)$,
there exists a $B>0$ such that for any $|b|\geq B$,
\[
f\left[b\right](z):=f(z)+\frac{1}{b}\int_{0}^{z}t e^{Q(t)}dt \in
 SF_{p,q}(\lambda) 
\]
is quasiconformally conjugate to $f$. For any $b\in\bC$, we write
\[
 F_b(z):=\frac{1}{b}f\left[ b\right](bz)=\lambda z(1+\frac{z}{2\lambda})
 +\frac{1}{b}h(bz),
\]
where $h$ is an entire function with $h(0)=h'(0)=0$.
\begin{proposition}\label{prop:hartogs}
If $f$ is linearizable at the origin, then
\[
 F_0(z)=\lambda z(1+\frac{z}{2\lambda})
\]
is also linearizable at the origin.
\end{proposition}

\begin{proof}
 We note that for $|b|\ge B$, $F_b$ is linearizable at the origin.
 As in the case of rational maps (cf. \cite{MSS} or \cite{McMullen:renorm}),
 we can show that the quasiconformal stability implies the $J$-stability.
 Hence in fact there exists an $M\ge 0$ such that for $B\le |b|\le 2B$,
 the Siegel disk of $F_b$ at the origin contains 
 $\{|z|\le M\}$. 
 Then the proposition follows by the same argument as in \cite{PMar},
 which is also explained in \cite{1-subhyp}. For completeness,
 we write the proof.

 Suppose that $J(F_0)$ intersects $\{|z|<M\}$.
 Then there exists a $z_1\in\bC$ with $0<|z_1|<M$ and $n>0$ such that
 $F_0^n(z_1)=z_1$ since $J(F_0)$ is the closure of the set of all
 repelling periodic points of $F_0$, which is true not only for
 rational functions but also entire functions. We set:
\[
 H(b,z):=\frac{z}{F_b^n(z)-z}:
 \{|b|<2B\}\times\{|z|<M\}\to \hbC,
\]
which depends meromorphically on each variables
and is uniformly continuous on $\{|b|\le 2B\}\times\{|z|\le M\}$.

 For $B<|b|<2B$,
 since $\{|z|\le M\}$ is contained in the Siegel disk of $F_b$ at the
 origin, $F_b$ has no periodic point there.
 Hence $H(b,z)$ is holomorphic on $\{B<|b|<2B\}\times\{|z|<M\}$.
 On the other hand, since $H(b,0)=1/(\lambda^{n}-1)$ is the constant
 independent of $|b|\le 2B$, there exists $0<m<M$ such that
 $H(b,z)$ is also holomorphic on $\{|b|<2B\}\times\{|z|<m\}$.

 By the Hartogs continuation theorem, $H(b,z)$ is actually holomorphic on
 $\{|b|<2B\}\times\{|z|<M\}$. This contradicts
 the assumption $F_0^n(z_1)=z_1$ and $0<|z_1|<M$.
\end{proof}

Hence $F_0$ is linearizable at the origin.
It follows from this and Theorem \ref{th:quad} that
$\alpha$ satisfies the Brjuno condition. \qed
\end{pfcritical}

\begin{pfsingularity}
We assume that $f(z)=\lambda\int_0^ze^{Q(t)}dt$, $Q(0)=0$.
Since $f$ is quasiconformally stable in $SF_{0,q}(\lambda)$,
there exists $B>0$ such that for any $|b|\geq B$,
\[
f\left[b\right](z):=\lambda\int_{0}^{z}e^{Q(t)+t/b}dt\in SF_{0,q}(\lambda)
\]
is quasiconformally conjugate to $f$. For any $b\in\bC$, we have
\[
 F_b(z):=\frac{1}{b}f\left[ b\right](bz)=
\lambda\int_0^ze^tdt+h(b,z),
\]
where $h$ is a holomorphic function on $\bC\times\bC$ with $h(0,\cdot)=0$.

\begin{proposition}
If $f$ is linearizable at the origin, then
\[
 F_0(z)=\lambda\int_0^ze^tdt
\]
is also linearizable at the origin.
\end{proposition}
 This can be proved by the same argument as in Proposition
 \ref{prop:hartogs}. \qed
\end{pfsingularity}

\section{Proof of the generalized Ma\~n\'e theorem}\label{sec:mane}
In this section, we show the generalized Ma\~n\'e theorem 
for every entire function with only finitely many critical points
and asymptotic values. 
In \cite{Mane93}, he showed it for rational functions (see also \cite{ST99}).
Throughout this section,
let $f$ be an entire function with only finitely many critical points
and asymptotic values.

\begin{theorem}\label{th:point}
 Let $M_f$ be the set of all asymptotic values and 
 recurrent critical points of $f$, and put $d(M_f):=\bigcup_{s\in M_f}d(s)$.
 
 Then there exists an $N\in\bN$ such that for every $x\in J(f)-d(M_f)$
 which is not a parabolic periodic point and for every $\epsilon>0$,
 there exists a connected neighborhood $U$ of $x$ such that for every $n\ge 0$
 and every connected component $V'$ of $f^{-n}(U)$ bounded in $\bC$,
\begin{enumerate}
\def\labelenumi{(\roman{enumi})}
 \item the spherical diameter of $V'$ is less than $\epsilon$ and
       $\deg(f^n:V'\to U)\le N$.
 \item For every $\epsilon_1>0$, there exists an $n_0\in\bN$
       such that for every $n>n_0$,
       the spherical diameter of $V'$ is less than $\epsilon_1$.
\end{enumerate}
\end{theorem}

Theorem \ref{th:point} can be shown by completely the same way
as Theorem 1.1 in \cite{ST99}.
The only difference is that we should exclude
not only the $\omega$-limit set of a recurrent critical point
but also the derived set of an asymptotic value.

Theorem \ref{th:invariant}, which is not needed to show Theorem
\ref{th:correspond}, follows from Theorem \ref{th:point} and
is also proved by almost the same argument as Theorem 1.2 in
\cite{ST99} although 
some extra argument for excluding unbounded iterated
preimages of the $U$ in Theorem \ref{th:invariant} is needed.

\begin{theorem}\label{th:invariant}
 Let $\Lambda\subset J(f)$ be compact
 and forward invariant, i.e., $f(\Lambda)\subset\Lambda$,
 and contain none of
 critical points, parabolic periodic points and asymptotic values.
 If $\Lambda\cap d(M_f)=\emptyset$, then it is expanding$;$
 i.e., there exists an $n_1>0$ such that for every $n\ge n_1$,
 $\min_{z\in\Lambda}|(f^n)'(z)|>1$.
\end{theorem}

\begin{proof}
 Assume that $\Lambda$ is not expanding. Then there exist
 $n_k\to\infty$, $\Lambda\ni z_k$, and $x\in\Lambda$ such that
 $|(f^{n_k})'(z_k)|\le 1$ and $\lim_{k\to\infty}f^{n_k}(z_k)=x$.
 If $x\in\bC-d(M_f)$,
 then $x$ satisfies the condition of Theorem \ref{th:point}.
 For every $\epsilon>0$ such that the spherical $\epsilon$-neighborhood
 $\Lambda_{\epsilon}$ of $\Lambda$ contains no critical point,
 let $U$ be a neighborhood of $x$ associated to $\epsilon$
 given by Theorem \ref{th:point}. Since
 $\lim_{k\to\infty}f^{n_k}(z_k)=x$, there exists $K(U)>0$ such that
 $f^{n_k}(z_k)\in U$ for every $k>K(U)$. 
 
 We show that there exists a $k>K(U)$ such that
 the component $V_k$ of $f^{-n_k}(U)$ containing $z_k$ is unbounded:
 otherwise, $f^{n_k}:V_k\to U$ is bijective for every $k>K(U)$
 since for $0\le j\le n_k$, $f^j(V_k)$ is a bounded connected component
 of $f^{-n_k+j}(U)$ intersecting $\Lambda$ 
 so is contained in $\Lambda_{\epsilon}$ by (i) in Theorem \ref{th:point}.
 Furthermore, every limit function $\phi$ on $U$ 
 of the single-valued branches $\{f^{-n_k}:U\to V_k\}_{k>K(U)}$ is a constant
 one by (ii) in Theorem \ref{th:point}. This contradicts that
 $|\phi'(x)|=\lim_{k\to\infty}|(f^{n_k})'(z_k)|^{-1}\ge 1$. 

 For an unbounded $V_k$, 
 there exists a $1\le j\le n_k$ such that $f^j(V_k)$ is bounded
 but $f^{j-1}(V_k)$ is unbounded. Then $f^j(V_k)$ contains an
 asymptotic value, and is contained in $\Lambda_{\epsilon}$.
 Consequently, if $\epsilon>0$ is small enough, 
 $\Lambda_{\epsilon}$ contains some asymptotic value.
 Since $f$ has only finitely many asymptotic values, 
 it implies $\Lambda$ itself contains some asymptotic value.
 This is a contradiction.
\end{proof}


Now we shall prove Theorem \ref{th:correspond}, the generalized Ma\~n\'e
theorem, stated in Section \ref{sec:subhyp}. For Cremer cycles,
we need some careful argument for finding bounded iterated preimages of the $U$
in Theorem \ref{th:invariant}.

\begin{proof}
 As in Section \ref{sec:subhyp},
 let $C$ be an irrationally indifferent cycle of period $p$,
 and $\Gamma=\Gamma_C$
 the boundary of the cycle of the Siegel disks associated with $C$
 if $C$ is Siegel, and the cycle $C$ itself otherwise.

 {\bfseries (Cremer Case)} Assume that $\Gamma$ is a Cremer cycle. Theorem
 \ref{th:invariant} implies that if $\Gamma$ contains no asymptotic
 value, then $\Gamma\subset d(M_f)$.
 We shall show that the latter always occurs.

 Assume that $\Gamma\cap d(M_f)=\emptyset$. Let $W$ be 
 the bounded and open spherical $\epsilon$-neighborhood of $\Gamma$.
 Since $\Gamma$ is a finite set and $f$ has only finitely many
 critical or asymptotic values,
 there exists an $\epsilon>0$ such that each connected component of
 $W-\Gamma$ is a spherical once-punctured disk whose puncture is in
 $\Gamma$ and $W-\Gamma$ contains none of critical or asymptotic values.
 Then $f:f^{-1}(W-\Gamma)\to W-\Gamma$ is unbranched.
 Hence $f$ maps each connected component $Y$ of $f^{-1}(W-\Gamma)$ onto
 a connected component $X$ of $W-\Gamma$ as a covering map,
 which is known to be isomorphic to that onto $\bD^*$
 given by either the logarithm or the $n$th root for some $n\in\bN$
 (cf. \cite{Forster}, Theorem 5.10). 
 $f:Y\to X$ is a logarithmic covering if and only if $Y$ is simply
 connected. We note that $f^{-1}(\Gamma)$ is the set of all punctures of
 $f^{-1}(W-\Gamma)$, and $\Gamma(\subset f^{-1}(\Gamma))$
 contains no critical point. Hence if $Y$ has its puncture
 in $\Gamma$, then $f:Y\to X$ is conformal,
 and $f$ gives a homeomorphism between the closures of $X$ and $Y$.
 In particular, since $X$ is bounded, so is $Y$.

 By filling the punctures, we conclude that
 every connected component of $f^{-1}(W)$
 that intersects $\Gamma$ is bounded and $f$ maps it onto $W$
 conformally. Let $W_1$ be the union of all these (only finitely many)
 connected components of $f^{-1}(W)$ that intersect $\Gamma$.

 For every $x\in\Gamma$, which satisfies the assumption of Theorem
 \ref{th:point}, there exists an open neighborhood $U_x$ of $\Gamma$
 associated to $\epsilon$ given by Theorem \ref{th:point}. Let
 $U:=\bigcup_{x\in\Gamma}U_x$. For $k\ge 0$, 
 let $V_k$ be the union of such components of $f^{-k}(U)$ that intersect
 $\Gamma$. By induction, we show $V_k\subset W$: for $k=0$, it is trivial.
 Assume that it is true for $k$. Since $V_{k+1}$ is the union of such
 components of $f^{-1}(V_k)$ that intersect $\Gamma$, 
 $V_{k+1}\subset W_1$, which is bounded, so
 $V_{k+1}\subset W$ by (i) in Theorem \ref{th:point}.

 Since $W$ contains no critical point,
 $f^k$ maps $V_k$ onto $U$ conformally. By (ii) in Theorem \ref{th:point},
 every limit function $\phi_x$ on $U_x$ of the single-valued branches
 $\{f^{-kp}:U_x\to V_{kp}\}_{k\ge 0}$ is constant.
 This contradicts that $|\phi_x'(x)|=\lim_{k\to\infty}|(f^{kp})'(x)|^{-1}=1$. 

 Hence $\Gamma\subset d(M_f)$, which concludes that $\Gamma\subset d(s)$ for
 some $s\in M_f$. 
 
 {\bfseries (Siegel case)} Next, assume that
 $\Gamma$ is the boundary of a cycle of the Siegel disks.

 Let $A(f)$ be the set of all asymptotic values.
 Assume that there exists $x\in(\Gamma-\bigcup_{n\ge 0}f^n(A(f)))-d(M_f)$.
 Then if a neighborhood of $x$ is small enough, it does not intersect
 $\bigcup_{n\ge 0}f^n(A(f))$. 
 Since $x$ satisfies the assumption of Theorem \ref{th:point},
 there exists an open neighborhood
 $U\subset\bC-\bigcup_{n\ge 0}f^n(A(f))$ of $x$ associated to some
 $\epsilon>0$ given by Theorem \ref{th:point}.
 Since $U$ intersects a Siegel disk,
 there exist $\epsilon_1>0$ and $n_k\to\infty$ such that for every $k\in\bN$,
 the spherical diameter of the connected component $V_k$
 of $f^{-n_k}(U)$ intersecting $\Gamma$ is more than $\epsilon_1$.
 On the other hand,
 since $U\subset\bC-\bigcup_{n\ge 0}f^n(A(f))$, $V_k$ is bounded for
 every $k\in\bN$. 
 Hence by (ii) in Theorem \ref{th:point}, there exists $n_0\in\bN$ such
 that the spherical diameter of $V_k$ is less than $\epsilon_1$ for every
 $k>n_0$. This is a contradiction.

 Hence $\Gamma-\bigcup_{n\ge 0}f^n(A(f))\subset d(M_f)$.
 Since the left hand side is dense in $\Gamma$ and $d(M_f)$ is closed, 
 it follows that $\Gamma\subset d(M_f)$.

 The proof of that $\Gamma\subset d(s)$ for some $s\in M_f$
 is almost the same as the proof of the original Ma\~n\'e theorem,
 so we give its outline and omit the details.
 
 Let $D$ be {\itshape one of} the Siegel disks,
 $\hat{\Gamma}$ the boundary of $D$ considered in $\hat{\bC}$,
 and $\mu$ the harmonic measure on $\hat{\Gamma}$
 with respect to the Siegel periodic point $z_0\in D$.
 The support of $\mu$ equals $\hat{\Gamma}$. In particular, 
 $\mu(\hat{\Gamma}-\{\infty\})>0$.

 The dynamical system
 $(\hat{\Gamma}-\{\infty\},
 f^p|(\hat{\Gamma}-\{\infty\}),\mu|(\hat{\Gamma}-\{\infty\}))$
 has such an ergodic property that
 if $f^p(B)\subset B$ for a Borel subset $B\subset\hat{\Gamma}-\{\infty\}$, 
 then $\mu(B)$ equals 0 or 1. In particular, for every $s\in\hbC$,
 $\mu((\hat{\Gamma}-\{\infty\})\cap d(s))$ equals 0 or 1.
 
 Consequently, since
 $\hat{\Gamma}-\{\infty\}\subset d(M_f)=\bigcup_{s\in M_f}d(s)$,
 there exists an $s\in M_f$ such that
 $\mu((\hat{\Gamma}-\{\infty\})\cap d(s))=1$. Hence
 $\hat{\Gamma}\cap d(s)$ is a $\mu$-full measure closed set in $\hbC$,
 so contains $\supp\mu=\hat{\Gamma}\supset\partial D$, which is the
 boundary of $D$ in $\bC$.
 Hence $d(s)\supset\partial D$
 so $d(s)\supset\overline{\bigcup_{n=1}^pf^n(\partial D)}=\Gamma$.
\end{proof}
  
\def\cprime{$'$}

\end{document}